\documentclass[11pt,a4paper]{article}
\usepackage[T1]{fontenc}
\usepackage[utf8]{inputenc}
\usepackage{lmodern}
\usepackage[margin=27mm]{geometry}
\usepackage{amsmath,amssymb,amsthm,mathtools}
\usepackage{booktabs,array}
\usepackage{tikz}
\usepackage{indentfirst}
\usepackage[colorlinks=true,linkcolor=blue,citecolor=blue,urlcolor=blue]{hyperref}
\allowdisplaybreaks[2]

\newtheorem{theorem}{Theorem}[section]
\newtheorem{proposition}[theorem]{Proposition}
\newtheorem{lemma}[theorem]{Lemma}
\newtheorem{corollary}[theorem]{Corollary}
\theoremstyle{definition}

\theoremstyle{remark}
\newtheorem{remark}[theorem]{Remark}

\newcommand{\HH}{\mathrm H}
\newcommand{\DD}{\mathrm D}
\newcommand{\ZZ}{\mathbb Z}
\newcommand{\QQ}{\mathbb Q}
\newcommand{\Dc}{D^b_c(R)}
\newcommand{\dc}{%
  \mathrel{\mathchoice
    {\dcl{0.16}}
    {\dcl{0.16}}
    {\dcl{0.11}}
    {\dcl{0.08}}
  }%
}
\newcommand{\dcl}[1]{\,%
  \mathrel{\raisebox{-0.15ex}{\hbox{%
  \begin{tikzpicture}[scale=#1, line width=0.6pt]
    \draw (0,0) -- (2,0) -- (1.2,1.75) -- (0.8,1.75) -- cycle;
    \draw (1.6,0)-- (0.8,1.75);
  \end{tikzpicture}}}%
}}
\newcommand{\Hd}{\widetilde{\HH}}
\newcommand{\coeur}{\operatorname{\textnormal{c\oe ur}}}
\newcommand{\nil}{\mathrm{nil}}
\newcommand{\SSpart}{\mathrm{ss}}
\newcommand{\tors}{\mathrm{tors}}
\newcommand{\ud}{\vee_{\mathrm u}}
\newcommand{\cd}{\vee_{\mathrm D}}
\newcommand{\NS}{\operatorname{NS}}
\newcommand{\Pic}{\operatorname{Pic}}
\newcommand{\Br}{\operatorname{Br}}
\newcommand{\DM}{\operatorname{DM}}
\newcommand{\Gal}{\operatorname{Gal}}
\newcommand{\Hom}{\operatorname{Hom}}
\newcommand{\Tot}{\operatorname{Tot}}
\newcommand{\TC}{\operatorname{TC}}
\newcommand{\im}{\operatorname{im}}

\title{Finite torsion in Hodge--Witt cohomology groups}
\author{Yuan Yang\\[2pt]
  \small BICMR,
  Peking University, Beijing, China\\
  \small\href{mailto:yy_@bicmr.pku.edu.cn}{\texttt{yy\_@bicmr.pku.edu.cn}}}
\date{}

\begin{document}
\maketitle
\vspace{-2em}

\begin{abstract}
For bounded coherent complexes over the Raynaud ring, we establish
natural Dieu\-donn\'e isomorphisms between Ekedahl's biduality pieces
$A_f^2$ and $A_{\SSpart}^1$, the nilpotent and semi-simple torsion of
the Illusie--Raynaud c{\oe}ur, and the corresponding parts of the finite
torsion $t^1$ of Ekedahl's diagonal cohomology. In bidegree $(0,i+1)$,
this torsion measures the non-smoothness of the $i$th flat Artin--Mazur
formal group; hence all such groups are smooth for Mazur--Ogus
varieties. For a smooth projective surface $S$, we describe the three
possible dual pairs via the Galois module $\NS(\bar{S})[p^\infty]$ and
the local-local and local-\'etale parts of the finite Picard quotient.
\end{abstract}

\section{The three descriptions and the comparison theorem}\label{sec:intro}

Let $k$ be a perfect field of characteristic $p>0$, let $W=W(k)$ and
$K=W[1/p]$, and let $\sigma$ be the Frobenius of $W$; we also write
$a^\sigma=\sigma(a)$. As in \cite[I, 1.1]{IR}, $R=R(k)=R^0\oplus R^1$
denotes the Raynaud ring: the graded $W$-algebra generated by $F$ and
$V$ in degree $0$ and by $d$ in degree $1$, subject to the relations
\[
 Fa=a^\sigma F,\quad aV=Va^\sigma,\quad FV=VF=p,\quad
 da=ad,\quad d^2=0,\quad FdV=d
\]
for $a\in W$. Its degree-zero part $R^0$ is the Dieudonn\'e ring, and by a
Dieudonn\'e module we mean a left $R^0$-module that is finitely generated
over $W$. We write $\Dc$ for the bounded derived category of graded left
$R$-modules with coherent cohomology
\cite[I, \S3]{IR}, \cite{EkI}. Coherence allows dominos; it does not
mean that the graded pieces are finitely generated over $W$.

For $M\in\Dc$, $\HH^j(M)^i$ denotes the part of internal (Raynaud)
degree $i$ of the $j$th cohomology module. Following
\cite[\S0]{EkI}, $M(a)[b]$ is the shift by $a$ in the $R$-grading and by
$b$ in the cohomological grading:
\[
 \HH^n\bigl(M(a)[b]\bigr)^m=\HH^{n+b}(M)^{m+a}.
\]
Thus $(-i)$ moves a module concentrated in internal degree $0$ to
internal degree $i$. ``Finite torsion'' always means finite
\emph{$W$-length}, not finiteness of the underlying set. The letter $T$
with two indices is reserved for the domino numbers $T^{i,j}$ of
Illusie--Raynaud, which appear only in Section~\ref{sec:surfaces};
Ekedahl's biduality filtration $T^a(L)$ in \S\ref{subsec:ekedahl}
carries a single index.

\subsection{The c\oe ur and its finite torsion}\label{subsec:coeur}

Let $L$ be a coherent graded left $R$-module. Following \cite[I, 1.4]{IR},
put $Z^i=\ker(d\colon L^i\to L^{i+1})$, $B^i=dL^{i-1}$,
\begin{align}
 V^{-\infty}Z^i&=\bigcap_{r\geq 0}V^{-r}Z^i,\qquad
 V^{-r}Z^i=\{x\in L^i : V^rx\in Z^i\},
 \label{eq:stable-cycles}\\
 F^{\infty}B^i&=\bigcup_{s\geq 0}F^sB^i,
 \label{eq:saturated-boundaries}
\end{align}
and define the c\oe ur of $L$ in degree $i$ by
\begin{equation}\label{eq:heart}
 \coeur(L^i)=V^{-\infty}Z^i/F^{\infty}B^i.
\end{equation}
When the module has to be specified we write $V^{-\infty}Z^i(L)$ and
$F^\infty B^i(L)$. The union in \eqref{eq:saturated-boundaries} is
increasing because $d=FdV$. Both $V^{-\infty}Z^i$ and $F^\infty B^i$
are $R^0$-submodules of $L^i$, and
$B^i\subset F^\infty B^i\subset V^{-\infty}Z^i\subset Z^i$
\cite[I, 1.4]{IR}. By the structure theory of coherent modules,
$\coeur(L^i)$ is a finitely generated $W$-module on which $V$ is
topologically nilpotent \cite[I, Theorem 2.9]{IR}; see also
\cite[Lemma 0.8]{EkI}.

Every Dieudonn\'e module $D$ of finite length has a canonical Fitting
decomposition \cite[I, (1.5.3)]{IR}
\begin{equation}\label{eq:fitting}
 D=D_{\nil}\oplus D_{\SSpart},\qquad
 D_{\nil}=\ker F^r,\quad D_{\SSpart}=\im F^r\quad(r\gg0).
\end{equation}
Here $\im F^r$ is a $W$-submodule because $k$ is perfect; $F$ is
nilpotent on $D_{\nil}$ and bijective on $D_{\SSpart}$; both summands
are $R^0$-submodules because $F$ and $V$ commute; and the decomposition
is functorial. As in the classical theory of $p$-linear maps
\cite[1.0.6--1.0.9]{KatzSGA7}, we call $D_{\nil}$ and $D_{\SSpart}$ the
\emph{nilpotent} and \emph{semi-simple} parts of $D$, and we call $D$ nilpotent or
semi-simple if it equals the corresponding part. ``Semi-simple'' refers
only to the bijectivity of $F$, not to semisimplicity as a
representation.

For $M\in\Dc$ and $i,j\in\ZZ$ put
\begin{equation}\label{eq:finite-heart}
 \coeur^{i,j}(M)=\coeur\bigl(\HH^j(M)^i\bigr),
\end{equation}
and let $\coeur^{i,j}(M)_{\tors}$ be its $p$-power torsion submodule.
This is a Dieudonn\'e module of finite length, and we write
\[
 \coeur^{i,j}(M)_{\tors}
 =\coeur^{i,j}(M)_{\nil}\oplus\coeur^{i,j}(M)_{\SSpart}
\]
for its decomposition \eqref{eq:fitting}. For a single coherent module
$L$ we use the same notation $\coeur(L^i)_{\tors}$,
$\coeur(L^i)_{\nil}$ and $\coeur(L^i)_{\SSpart}$.

\subsection{Ekedahl's biduality filtration}\label{subsec:ekedahl}

Let $\DD$ be Ekedahl's duality functor on $\Dc$, and put
$\DD^a(L)=\HH^a(\DD(L))$ \cite[IV, 2.1]{EkI}. For a coherent graded
left $R$-module $L$, the biduality spectral sequence gives a filtration
\[
 0\subset T^2(L)\subset T^1(L)\subset T^0(L)=L,
 \qquad A^a(L)=T^a(L)/T^{a+1}(L),
\]
with $A^a(L)\simeq\DD^a\DD^a(L)$ \cite[IV, 2.3--2.4]{EkI}. The
filtration has length three, so $A^2(L)=T^2(L)$. As in
\cite[IV, 3.2--3.3]{EkI}, the two finite pieces are
\begin{equation}\label{eq:A-pieces}
 A_f^2(L)^i=\coeur\bigl(A^2(L)^i\bigr),\qquad
 A_{\SSpart}^1(L)=A^1(L)_{\tors}.
\end{equation}
Both have finite $W$-length; $F$ is nilpotent on $A^2_f(L)$ and
bijective on $A^1_{\SSpart}(L)$ \cite[IV, Theorem 3.3]{EkI}. Moreover,
\begin{equation}\label{eq:duality-quotient}
 T^2(L)\subset L_{\tors},\qquad
 L_{\tors}/T^2(L)=A^1(L)_{\tors}=A_{\SSpart}^1(L),
\end{equation}
where $L_{\tors}$ is the $p$-power torsion of $L$. Indeed, $T^2(L)$ is
$p$-power torsion \cite[IV, Proposition 2.4 and Lemma 2.4.1]{EkI}, and
the filtration gives an exact sequence
\[
 0\longrightarrow A^1(L)\longrightarrow L/T^2(L)
 \longrightarrow A^0(L)\longrightarrow 0
\]
with $A^0(L)$ torsion-free \cite[IV, Theorem 3.3]{EkI}. Hence the
torsion of $L/T^2(L)$, which is $L_{\tors}/T^2(L)$, equals
$A^1(L)_{\tors}$. The differential vanishes on $A^1_{\SSpart}(L)$: $V$
is nilpotent there, so $d=F^rdV^r$ is zero for $r\gg0$.

\subsection{Diagonal finite torsion}\label{subsec:diagonal}

Write $\dc$ for the heart of Ekedahl's diagonal $t$-structure on $\Dc$,
and $\Hd^n$ for its cohomology functors \cite{EkDiag}.
The ordinary and diagonal $t$-structures commute. If $A\in\dc$, its
ordinary cohomology satisfies
\begin{equation}\label{eq:diag-support}
 \HH^q(A)^a=0\quad\text{unless }a\in\{-q,1-q\},\qquad
 F^\infty B^{1-q}\bigl(\HH^q(A)\bigr)=\HH^q(A)^{1-q}.
\end{equation}
Thus $\HH^q(A)(-q)$ is a normalized diagonal module: it is concentrated
in degrees $0$ and $1$, and its degree-one term consists of saturated
boundaries.

Ekedahl's diagonal duality gives endofunctors $t^1$ and $t^0$ of $\dc$
and, for every object $A$ of $\dc$, a natural exact sequence
\begin{equation}\label{eq:finite-diag-sequence}
 0\longrightarrow t^1A\longrightarrow A\longrightarrow t^0A
 \longrightarrow 0,
\end{equation}
in which $t^1A$ is the largest finite-torsion subobject of $A$ and
$t^0A$ has no finite torsion; in terms of the diagonal duality,
$t^1=\widetilde{\DD}^{1}\widetilde{\DD}^{1}$ and
$t^0=\widetilde{\DD}^{0}\widetilde{\DD}^{0}$
\cite[I, Proposition 1.11]{EkDiag}. The superscripts of $t^1$ and $t^0$
are part of the notation; they are not cohomological degrees.

For $M\in\Dc$ and $i,j\in\ZZ$, the finite torsion of diagonal
cohomology in bidegree $(i,j)$ is the finite-length Dieudonn\'e module
(Lemma~\ref{lem:diagonal-coeur})
\begin{equation}\label{eq:diagonal-finite-bidegree}
 \HH^{-i}\bigl(t^1\Hd^{i+j}(M)\bigr)^i,
\end{equation}
with nilpotent and semi-simple parts
$\HH^{-i}\bigl(t^1\Hd^{i+j}(M)\bigr)^i_{\nil}$ and
$\HH^{-i}\bigl(t^1\Hd^{i+j}(M)\bigr)^i_{\SSpart}$ in the sense of
\eqref{eq:fitting}.

\begin{theorem}[Comparison of finite torsion]\label{thm:comparison}
Let $M\in\Dc$ and $i,j\in\ZZ$. There are natural isomorphisms of
finite-length Dieudonn\'e modules
\begin{align}
 A_f^2\bigl(\HH^j(M)\bigr)^i
 &\simeq \coeur^{i,j}(M)_{\nil}
 \simeq \HH^{-i}\bigl(t^1\Hd^{i+j}(M)\bigr)^i_{\nil},
 \label{eq:main-nil}\\
 A_{\SSpart}^1\bigl(\HH^j(M)\bigr)^i
 &\simeq \coeur^{i,j}(M)_{\SSpart}
 \simeq \HH^{-i}\bigl(t^1\Hd^{i+j}(M)\bigr)^i_{\SSpart}.
 \label{eq:main-ss}
\end{align}
In particular, these isomorphisms commute with $F$ and $V$; they are not
merely equalities of lengths or isomorphisms of $W$-modules.
\end{theorem}

By Theorem~\ref{thm:comparison}, if one of the three descriptions of
finite torsion vanishes in a given bidegree $(i,j)$, then so do the
other two. When this happens, that is, when
$\coeur^{i,j}(M)_{\tors}=0$, we say that $M$ has \emph{no finite
torsion in bidegree $(i,j)$}.

Section~\ref{sec:ordinary} proves the comparison with Ekedahl's pieces
directly from the structure of coherent modules and Ekedahl's biduality
filtration, and Section~\ref{sec:diagonal} deduces the diagonal
comparison from Ekedahl's finite-torsion criterion, keeping track of
both degrees. Section~\ref{sec:artin-mazur} interprets the finite
torsion in internal degree zero as the obstruction to smoothness of the
flat Artin--Mazur formal groups, deduces the corresponding
prorepresentability criterion for the \'etale Artin--Mazur functors, in
particular for the formal Brauer group, and applies this to
Mazur--Ogus varieties. Section~\ref{sec:surfaces} specializes to
surfaces: the Kummer sequence and logarithmic de Rham--Witt cohomology
identify the N\'eron--Severi contribution, the Cartier-module
comparison of Section~\ref{sec:artin-mazur} identifies the connected
Picard contributions, and geometric duality determines their partners.
The two duality conventions are kept apart.

\section{Comparison with Ekedahl's pieces}\label{sec:ordinary}

Fix a coherent graded left $R$-module $L$, and put
\[
 P=L_{\tors},\qquad Q=T^2(L),\qquad\text{so that }Q\subset P\subset L
\]
by \eqref{eq:duality-quotient}. By the structure theory of coherent
modules \cite[I, Theorem 2.9 and \S3]{IR} (see also \cite[IV, p.~223]{EkI}
for the torsion and torsion-free quotients), there is an integer $e\geq0$
with $p^eP=0$, and $L/P$ is a finitely generated torsion-free $W$-module
with zero differential on which $V$ is topologically nilpotent. In
particular $dL\subset P$ and $p^e\,dL=0$. By
\eqref{eq:duality-quotient}, $P/Q=A_{\SSpart}^1(L)$ has finite length
and bijective $F$, so $V=pF^{-1}$ is nilpotent on it.

We also use the domino of $L$ in degrees $(i,i+1)$ \cite[I, 2.D]{IR}:
\[
 L^i/V^{-\infty}Z^i\ \xrightarrow{\ d\ }\ F^\infty B^{i+1}.
\]
$F$ is nilpotent on its degree-$i$ term $L^i/V^{-\infty}Z^i$. Indeed,
$FdV=d$ gives $dF=pFd$,
hence $dF^m=p^mF^md$, and for $x\in L^i$ and $n\geq e$,
\[
 dV^rF^nx=
 \begin{cases}
 p^nF^{n-r}dx=0, & r\leq n,\\
 p^n\,dV^{r-n}x=0, & r>n.
 \end{cases}
\]
Thus $F^nL^i\subset V^{-\infty}Z^i$ for $n\geq e$. The argument only
uses $p^e\,dL=0$, so it applies equally to every graded $R$-submodule
of $L$.

\begin{lemma}\label{lem:saturated-boundaries}
The inclusions $Q\subset P\subset L$ induce equalities of submodules
\[
 F^\infty B^i(Q)=F^\infty B^i(P)=F^\infty B^i(L)
\]
and injective maps
\[
 \coeur(Q^i)\lhook\joinrel\longrightarrow\coeur(P^i)
 \lhook\joinrel\longrightarrow\coeur(L^i).
\]
\end{lemma}

\begin{proof}
Let $x\in L^{i-1}$. Since $V$ is topologically nilpotent on $L/P$,
there is an integer $r$ with
\[
 V^rx=p^ey+u,\qquad y\in L^{i-1},\quad u\in P^{i-1}.
\]
Iterating $d=FdV$ and using $p^e\,dL=0$, we obtain
\[
 dx=F^rdV^rx=F^r(p^e\,dy+du)=F^rdu\in F^\infty B^i(P).
\]
Since $F^\infty B^i(P)$ is $F$-stable, this proves
$F^\infty B^i(L)\subset F^\infty B^i(P)$; the reverse inclusion is
clear. Next, $V$ is nilpotent on $P/Q$, so $V^rP\subset Q$ for some
$r$, and for $x\in P^{i-1}$ we get $dx=F^rdV^rx\in F^\infty B^i(Q)$.
This proves the other equality in the same way.

Finally, $V^{-\infty}Z^i(J)=J^i\cap V^{-\infty}Z^i(L)$ for every graded
$R$-submodule $J\subset L$. Dividing these inclusions of stable cycles
by their common submodule $F^\infty B^i$ gives the injections of
c\oe urs.
\end{proof}

\begin{proposition}\label{prop:ordinary-comparison}
There are natural identifications
\[
 \coeur(P^i)=\coeur(L^i)_{\tors},\qquad
 \coeur(Q^i)=\coeur(L^i)_{\nil},\qquad
 P^i/Q^i\simeq\coeur(L^i)_{\SSpart}.
\]
\end{proposition}

\begin{proof}
By Lemma~\ref{lem:saturated-boundaries},
\[
 \coeur(L^i)/\coeur(P^i)=V^{-\infty}Z^i(L)/V^{-\infty}Z^i(P).
\]
Since $V^{-\infty}Z^i(P)=V^{-\infty}Z^i(L)\cap P^i$, this quotient
embeds in $L^i/P^i$, which is torsion-free. On the other hand,
$\coeur(P^i)$ is a subquotient of $P^i$, hence torsion, and it is
finitely generated over $W$ as a submodule of $\coeur(L^i)$. It is
therefore exactly the torsion submodule of $\coeur(L^i)$.

Applying the same lemma to $Q\subset P$ gives
\begin{equation}\label{eq:coeur-quotient}
 \coeur(P^i)/\coeur(Q^i)=V^{-\infty}Z^i(P)/V^{-\infty}Z^i(Q).
\end{equation}
The inclusions of stable cycles and their quotients give an exact
sequence
\begin{equation}\label{eq:stable-cycle-exact}
 0\longrightarrow\frac{V^{-\infty}Z^i(P)}{V^{-\infty}Z^i(Q)}
 \longrightarrow\frac{P^i}{Q^i}
 \longrightarrow
 \frac{P^i/V^{-\infty}Z^i(P)}{Q^i/V^{-\infty}Z^i(Q)}
 \longrightarrow 0.
\end{equation}
The last term is a quotient of the degree-$i$ term
$P^i/V^{-\infty}Z^i(P)$ of the domino of $P$, so $F$ is nilpotent on it. It is also a quotient of the finite-length module
$P^i/Q^i$, on which $F$ is bijective; hence $F$ is surjective on it and,
by finite length and the perfectness of $k$, bijective. A module on which
$F$ is both nilpotent and bijective is zero. By \eqref{eq:A-pieces} and
\eqref{eq:duality-quotient}, $\coeur(Q^i)=A_f^2(L)^i$ and
$P^i/Q^i=A_{\SSpart}^1(L)^i$, so \eqref{eq:coeur-quotient} and
\eqref{eq:stable-cycle-exact} give an exact sequence
\begin{equation}\label{eq:nil-ss-sequence}
 0\longrightarrow A_f^2(L)^i\longrightarrow\coeur(L^i)_{\tors}
 \longrightarrow A_{\SSpart}^1(L)^i\longrightarrow 0.
\end{equation}
The first term has nilpotent $F$ and the last has bijective $F$, so the
canonical decomposition \eqref{eq:fitting} identifies them with the
nilpotent and semi-simple parts of $\coeur(L^i)_{\tors}$. Every map used
is induced by the natural filtrations and is $R^0$-linear.
\end{proof}

\begin{remark}
The c\oe ur is not an exact functor on coherent modules, so the
injectivity in Lemma~\ref{lem:saturated-boundaries} is not automatic:
it relies on the equality of saturated boundaries for
$Q\subset P\subset L$.
\end{remark}

\section{Comparison with diagonal finite torsion}\label{sec:diagonal}

We first determine what $t^1$ does to ordinary cohomology, using
Ekedahl's description of finite torsion
\cite[I, Theorem 1.12]{EkDiag}:
\begin{enumerate}
\item $A\in\dc$ is of finite torsion if and only if every $\HH^q(A)$
is a finite-length Dieudonn\'e module concentrated in internal degree
$-q$;
\item $A\in\dc$ has no finite torsion if and only if no normalized
diagonal module $\HH^q(A)(-q)$ has finite torsion, and a normalized
diagonal module has no finite torsion if and only if the c\oe ur of its
degree-zero term is torsion-free.
\end{enumerate}
These criteria concern finite torsion only; the $p$-torsion carried by
dominos is not detected by $t^1$.

\begin{lemma}\label{lem:diagonal-coeur}
Let $A\in\dc$ and $q\in\ZZ$. Then $\HH^q(t^1A)$ is concentrated in
internal degree $-q$, and there is a natural $R^0$-linear
identification
\begin{equation}\label{eq:diagonal-coeur}
 \HH^q(t^1A)^{-q}\simeq\coeur\bigl(\HH^q(A)^{-q}\bigr)_{\tors}.
\end{equation}
\end{lemma}

\begin{proof}
The first assertion is criterion (1). Apply ordinary cohomology to
\eqref{eq:finite-diag-sequence}. The map
$\HH^{q-1}(t^0A)\to\HH^q(t^1A)$ is zero, because the source lives in
internal degrees $1-q$ and $2-q$ and the target only in degree $-q$.
For the same reason $\HH^q(t^0A)\to\HH^{q+1}(t^1A)$ vanishes. Hence
there are short exact sequences of graded left $R$-modules
\begin{equation}\label{eq:ordinary-exact-from-diagonal}
 0\longrightarrow\HH^q(t^1A)\longrightarrow\HH^q(A)
 \longrightarrow\HH^q(t^0A)\longrightarrow 0.
\end{equation}
The last two modules live in internal degrees $-q$ and $1-q$. In
particular nothing arrives in degree $-q$ under $d$, so their c\oe urs
in that degree are their submodules $V^{-\infty}Z^{-q}$ of stable
cycles.

The first term of \eqref{eq:ordinary-exact-from-diagonal} has finite
length and is concentrated in degree $-q$, so its image lies in
$\coeur(\HH^q(A)^{-q})_{\tors}$. Conversely, a torsion stable cycle of
$\HH^q(A)^{-q}$ maps to a torsion element of
$\coeur(\HH^q(t^0A)^{-q})$, which is torsion-free by criterion (2);
so its image is zero. Exactness proves \eqref{eq:diagonal-coeur}. All
the maps involved commute with $F$ and $V$.
\end{proof}

We now relate the diagonal cohomology of a complex to its individual
Hodge--Witt positions. For a coherent graded left $R$-module $L$, the
radical filtration defining the diagonal $t$-structure is
\begin{equation}\label{eq:radical-filtration}
 (\widetilde\tau_{\leqslant i}L)^a=
 \begin{cases}
 L^a, & a\leq i,\\
 F^\infty B^{i+1}(L), & a=i+1,\\
 0, & a>i+1.
 \end{cases}
\end{equation}
Its $i$th layer
$S^iL=\widetilde\tau_{\leqslant i}L/\widetilde\tau_{\leqslant i-1}L$ is
the two-term graded left $R$-module
\begin{equation}\label{eq:radical-layer}
 (S^iL)^i=L^i/F^\infty B^i(L)
 \ \xrightarrow{\ d\ }\
 (S^iL)^{i+1}=F^\infty B^{i+1}(L).
\end{equation}
In particular,
\begin{equation}\label{eq:layer-coeur}
 \coeur\bigl((S^iL)^i\bigr)=\coeur(L^i).
\end{equation}
Indeed, the stable cycles in the degree-$i$ term of
\eqref{eq:radical-layer} are exactly $V^{-\infty}Z^i(L)/F^\infty B^i(L)$, and $S^iL$ is zero in
degree $i-1$.

Commutation of the ordinary and diagonal truncations gives the
canonical layer identity
\begin{equation}\label{eq:commuting-truncations}
 \HH^{-i}\bigl(\Hd^{i+j}(M)\bigr)\simeq S^i\bigl(\HH^j(M)\bigr).
\end{equation}
To check the indices, take $L$ in cohomological degree $0$: its $i$th
diagonal layer is $S^iL[i]$, whose only ordinary cohomology lies in
degree $-i$; a summand in cohomological degree $j$ therefore contributes
to diagonal degree $i+j$. The general identity follows from the
commuting truncations, without assuming that $M$ splits as a sum of its
cohomology modules.

\begin{proposition}\label{prop:diagonal-comparison}
For every $M\in\Dc$ and $i,j\in\ZZ$, there is a natural $R^0$-linear
isomorphism
\[
 \HH^{-i}\bigl(t^1\Hd^{i+j}(M)\bigr)^i\simeq\coeur^{i,j}(M)_{\tors}.
\]
\end{proposition}

\begin{proof}
Apply Lemma~\ref{lem:diagonal-coeur} to $A=\Hd^{i+j}(M)$ and $q=-i$,
then use \eqref{eq:commuting-truncations} and \eqref{eq:layer-coeur}:
\begin{align*}
 \HH^{-i}\bigl(t^1\Hd^{i+j}(M)\bigr)^i
 &\simeq\coeur\Bigl(\HH^{-i}\bigl(\Hd^{i+j}(M)\bigr)^i\Bigr)_{\tors}\\
 &\simeq\coeur\Bigl(S^i\bigl(\HH^j(M)\bigr)^i\Bigr)_{\tors}\\
 &\simeq\coeur\bigl(\HH^j(M)^i\bigr)_{\tors}.
\end{align*}
Taking nilpotent and semi-simple parts gives the comparison for each
part separately.
\end{proof}

Propositions~\ref{prop:ordinary-comparison} and
\ref{prop:diagonal-comparison} prove Theorem~\ref{thm:comparison}. They
also give the comparison maps explicitly; no length count is needed.

\begin{corollary}\label{cor:whole-diagonal}
For every $n$, the finite torsion of diagonal cohomology is
\begin{equation}\label{eq:whole-diagonal}
 t^1\Hd^n(M)\simeq\bigoplus_{i+j=n}\coeur^{i,j}(M)_{\tors}(-i)[i],
\end{equation}
where each $\coeur^{i,j}(M)_{\tors}$ on the right is first placed in
internal and cohomological degree zero. In particular, there is an
isomorphism of $W$-modules
\begin{equation}\label{eq:finite-total}
 \HH^0\Tot\bigl(t^1\Hd^n(M)\bigr)
 \simeq\bigoplus_{i+j=n}\coeur^{i,j}(M)_{\tors},
\end{equation}
where $\Tot$ is the associated simple complex.
\end{corollary}

\begin{proof}
Ekedahl's finite-torsion criterion also decomposes a finite-torsion
object of $\dc$ into its ordinary cohomology pieces
\cite[I, Theorem 1.12]{EkDiag}. Insert
Proposition~\ref{prop:diagonal-comparison}: a piece in cohomological
degree $-i$ has internal degree $i$, which gives precisely the shift
$(-i)[i]$. Each such piece has total degree zero, which proves
\eqref{eq:finite-total}.
\end{proof}

\begin{remark}\label{rem:frobenius}
The isomorphisms of Theorem~\ref{thm:comparison} preserve the Raynaud
operator $F$ in each fixed internal degree; no Frobenius twist enters
\eqref{eq:diagonal-finite-bidegree}. Geometrically, crystalline
Frobenius in internal degree $i$ is $p^iF$. This is why
\eqref{eq:finite-total} is stated only as an identity of $W$-modules:
it does not identify crystalline Frobenius with the unscaled $F$ on
every summand.
\end{remark}

\section{Finite torsion and Artin--Mazur formal groups}\label{sec:artin-mazur}

In this section $X$ is a smooth proper scheme over the perfect field
$k$. Put
\[
 M_X=\mathrm R\Gamma(X,W\Omega_X^\bullet)\in\Dc,\qquad
 \coeur^{a,b}(X)=\coeur^{a,b}(M_X),
\]
with torsion, nilpotent and semi-simple parts as in
\S\ref{subsec:coeur}. Thus $\coeur^{0,b}(X)$ is the c\oe ur of
$\HH^b(M_X)^0=\HH^b(X,W\mathcal O_X)$.

For an Artinian $k$-algebra $B$, the Artin--Mazur functor \cite{AM} is
\begin{equation}\label{eq:etale-artin-mazur}
 \Phi^i_{\mathrm{\acute et}}(X,\mathbf G_m)(B)
 =\ker\bigl(\HH^i_{\mathrm{\acute et}}(X_B,\mathbf G_m)
 \longrightarrow
 \HH^i_{\mathrm{\acute et}}(X_{B_{\mathrm{red}}},\mathbf G_m)\bigr).
\end{equation}
A functor of this kind is called \emph{prorepresentable} if it is
representable by a formal group. Let
\[
 G_i=\Phi^i_{\mathrm{fl}}(X,\mathbf G_m)
\]
be the fppf sheafification of $\Phi^i_{\mathrm{\acute et}}(X,\mathbf G_m)$
on the category of Artinian $k$-algebras. On each $X_B$, \'etale and
fppf cohomology with coefficients in $\mathbf G_m$ agree, so the two
functors differ only through this sheafification. The functor $G_i$ is
always prorepresentable, by a connected commutative formal group with
finite-dimensional tangent space \cite[Theorem 10.8]{BO}. In this
section, a formal group is called smooth if it is formally smooth over
$k$.

\subsection{The finite obstruction to smoothness}\label{subsec:obstruction}

Since $k$ is perfect, the reduced subgroup $(G_i)_{\mathrm{red}}$ is a
smooth formal group, and
\[
 Q_i=G_i/(G_i)_{\mathrm{red}}
\]
is a finite infinitesimal commutative group scheme; thus $G_i$ is smooth
if and only if $Q_i=0$. For instance, $G_0$ is smooth and $Q_0=0$: the
$k$-algebra $\HH^0(X,\mathcal O_X)$ is finite \'etale, and
$G_0=\Phi^0_{\mathrm{\acute et}}(X,\mathbf G_m)$ is the formal completion
at the identity of its group of units.

Let $\DM_k$ be the covariant Dieudonn\'e functor over $k$, namely the
contravariant Dieudonn\'e module of the Cartier dual. A finite
commutative group scheme of $p$-power order over $k$ is canonically the
direct sum of three parts, indicated by the subscripts $\mathrm{el}$,
$\mathrm{ll}$ and $\mathrm{le}$: the first letter records whether the
part is \'etale or local, and the second does the same for its Cartier
dual. No \'etale part with \'etale Cartier dual occurs in $p$-power
order, and the decomposition is compatible with extension of perfect
fields. Both $F$ and $V$ are nilpotent on $\DM_k$ of a local-local
group scheme; $F$ is bijective and $V$ nilpotent on $\DM_k$ of a
local-\'etale one; and $V$ is bijective and $F$ nilpotent on $\DM_k$ of
an \'etale-local one. Being infinitesimal, $Q_i$ is the sum
$(Q_i)_{\mathrm{ll}}\oplus(Q_i)_{\mathrm{le}}$ of its local-local and
local-\'etale parts, and $V$ is nilpotent on $\DM_k(Q_i)$.

Cartier theory identifies this Dieudonn\'e module with the first
derived functor of the Cartier module of typical curves $\TC$:
\begin{equation}\label{eq:derived-cartier-obstruction}
 R^1\TC(G_i)\simeq\DM_k(Q_i)
 \qquad\text{\cite[Remark 11.10]{BO}}.
\end{equation}
In particular, $G_i$ is smooth if and only if $R^1\TC(G_i)=0$.

\begin{theorem}[Smoothness and finite torsion]\label{thm:artin-mazur-finite-torsion}
For every $i\geq0$, there are natural isomorphisms of finite-length
Dieudonn\'e modules
\begin{equation}\label{eq:artin-mazur-finite-torsion}
 \coeur^{0,i+1}(X)_{\tors}=\HH^{i+1}(X,W\mathcal O_X)[V^\infty]
 \simeq R^1\TC(G_i)\simeq\DM_k(Q_i).
\end{equation}
Consequently, $\Phi^i_{\mathrm{fl}}(X,\mathbf G_m)$ is smooth if and
only if $M_X$ has no finite torsion in bidegree \mbox{$(0,i+1)$}. More
precisely,
\[
 \coeur^{0,i+1}(X)_{\nil}\simeq\DM_k\bigl((Q_i)_{\mathrm{ll}}\bigr),
 \qquad
 \coeur^{0,i+1}(X)_{\SSpart}\simeq\DM_k\bigl((Q_i)_{\mathrm{le}}\bigr).
\]
\end{theorem}

\begin{proof}
Ekedahl's flat Artin--Mazur comparison
\cite[III, Proposition 8.1]{EkII} gives a short exact sequence
\begin{equation}\label{eq:general-cartier-sequence}
 0\longrightarrow R^1\TC(G_i)\longrightarrow\HH^{i+1}(X,W\mathcal O_X)
 \longrightarrow\TC(G_{i+1})\longrightarrow 0.
\end{equation}
By \eqref{eq:derived-cartier-obstruction}, $V$ is nilpotent on the first
term, and it is injective on the Cartier module $\TC(G_{i+1})$. Hence
the first term is precisely the $V$-power torsion of the middle term.

We identify this $V$-power torsion with the c\oe ur torsion. Nothing
arrives in internal degree zero under $d$, so the c\oe ur in that degree
is the module of stable cycles. If $V^mx=0$, then $p^mx=F^mV^mx=0$ and,
for every $r\geq0$ and $n\geq m$,
\[
 dV^rx=F^ndV^{r+n}x=0;
\]
thus $x$ is a torsion stable cycle. Conversely, the torsion of the
c\oe ur has finite length and $V$ is topologically nilpotent on it, so
$V$ is nilpotent on it. This proves the equality in
\eqref{eq:artin-mazur-finite-torsion}, and the two isomorphisms follow
from \eqref{eq:general-cartier-sequence} and
\eqref{eq:derived-cartier-obstruction}.

The covariant Dieudonn\'e functor detects the zero group scheme, so
$\DM_k(Q_i)=0$ if and only if $Q_i=0$, that is, if and only if $G_i$ is
smooth; by \eqref{eq:artin-mazur-finite-torsion}, this happens if and
only if $\coeur^{0,i+1}(X)_{\tors}=0$. Finally, $F$ is
nilpotent on $\DM_k((Q_i)_{\mathrm{ll}})$ and bijective on
$\DM_k((Q_i)_{\mathrm{le}})$, which gives the last assertion.
\end{proof}

Theorem~\ref{thm:artin-mazur-finite-torsion} concerns only the finite
torsion of the c\oe ur: $\HH^{i+1}(X,W\mathcal O_X)$ may contain further
$p$-torsion coming from dominos. Since $G_0$ is smooth, the case $i=0$
gives
\begin{equation}\label{eq:H1-WO}
 \coeur^{0,1}(X)_{\tors}=\HH^1(X,W\mathcal O_X)[V^\infty]=0.
\end{equation}

\subsection{The \'etale Artin--Mazur functors}

The passage from the flat formal group $G_i$ to the \'etale functor
$\Phi^i_{\mathrm{\acute et}}(X,\mathbf G_m)$ involves the preceding
cohomological degree. The precise criterion is
\cite[Proposition 10.11]{BO}; the sufficient prorepresentability
criterion goes back to \cite[II, Proposition 1.8 and Corollary 2.11]{AM}.

\begin{corollary}\label{cor:etale-artin-mazur}
Let $i\geq1$. The functor $\Phi^i_{\mathrm{\acute et}}(X,\mathbf G_m)$ is
prorepresentable if and only if $\coeur^{0,i}(X)_{\tors}=0$, and it is
prorepresentable and smooth if and only if
\begin{equation}\label{eq:etale-two-degrees}
 \coeur^{0,i}(X)_{\tors}=\coeur^{0,i+1}(X)_{\tors}=0.
\end{equation}
Whenever it is prorepresentable, the natural map
\[
 \Phi^i_{\mathrm{\acute et}}(X,\mathbf G_m)
 \longrightarrow\Phi^i_{\mathrm{fl}}(X,\mathbf G_m)
\]
is an isomorphism. In particular, if all the flat Artin--Mazur formal
groups are smooth, then all the \'etale Artin--Mazur functors are
prorepresentable and smooth.
\end{corollary}

\begin{proof}
By \cite[Propositions 10.10--10.11]{BO},
$\Phi^i_{\mathrm{\acute et}}(X,\mathbf G_m)$ is prorepresentable if and
only if $G_{i-1}$ is smooth, and this is also equivalent to the
comparison map $\Phi^i_{\mathrm{\acute et}}(X,\mathbf G_m)\to G_i$ being
an isomorphism. Theorem~\ref{thm:artin-mazur-finite-torsion}, applied
to $i-1$ and to $i$, proves the first two assertions. The last one
follows because $G_0=\Phi^0_{\mathrm{\acute et}}(X,\mathbf G_m)$ is
smooth.
\end{proof}

Thus smoothness of $\Phi^i_{\mathrm{fl}}(X,\mathbf G_m)$ alone
guarantees prorepresentability of
$\Phi^{i+1}_{\mathrm{\acute et}}(X,\mathbf G_m)$, while
prorepresentability and smoothness of
$\Phi^i_{\mathrm{\acute et}}(X,\mathbf G_m)$ are governed by the two
adjacent finite-torsion modules in \eqref{eq:etale-two-degrees}. For
$i=1$ and $i=2$, $\Phi^i_{\mathrm{\acute et}}(X,\mathbf G_m)$ is the
formal Picard group and the formal Brauer group of \cite{AM}. By
\eqref{eq:H1-WO}, the formal Picard group is always prorepresentable,
and the formal Brauer group is prorepresentable if and only if the
formal Picard group is smooth, that is, if and only if
$\coeur^{0,2}(X)_{\tors}=0$.

\subsection{Mazur--Ogus varieties}

Recall that $X$ is a \emph{Mazur--Ogus} variety if all the crystalline
cohomology groups $\HH^n_{\mathrm{cris}}(X/W)$ are torsion-free and the
Hodge--de Rham spectral sequence of $X$ degenerates at $E_1$.

\begin{proposition}\label{prop:mazur-ogus-finite-torsion}
If $X$ is Mazur--Ogus, then
\[
 t^1\Hd^n(M_X)=0\quad\text{for all }n,\qquad
 \coeur^{a,b}(X)_{\tors}=0\quad\text{for all }a,b.
\]
Thus $M_X$ has no finite torsion in any bidegree.
\end{proposition}

\begin{proof}
Put $A_n=\Hd^n(M_X)$. By \cite[IV, Theorem 1.2]{EkDiag}, $A_n$ is a
Mazur--Ogus object of $\dc$, so by \cite[I, Lemma 4.4]{EkDiag} the
complex $\Tot(A_n)$ is concentrated in degree zero and
$\HH^0\Tot(A_n)$ is torsion-free over $W$. For every $A\in\dc$, the
module $\HH^q(A)^a$ contributes to $\Tot(A)$ in degree $q+a$, so
\eqref{eq:diag-support} shows that $\Tot(A)$ has cohomology only in
degrees $0$ and $1$; in particular $\HH^{-1}\Tot(t^0A_n)=0$. Applying
$\Tot$ to the sequence \eqref{eq:finite-diag-sequence} for $A_n$, we
therefore obtain an injection
\[
 \HH^0\Tot(t^1A_n)\lhook\joinrel\longrightarrow\HH^0\Tot(A_n).
\]
By \eqref{eq:finite-total}, its source is the finite-length $W$-module
$\bigoplus_{a+b=n}\coeur^{a,b}(X)_{\tors}$, whereas its target is
torsion-free. Hence $\coeur^{a,b}(X)_{\tors}=0$ whenever $a+b=n$, and
\eqref{eq:whole-diagonal} gives $t^1A_n=0$.
\end{proof}

\begin{corollary}\label{cor:mazur-ogus-artin-mazur}
Let $X$ be Mazur--Ogus. For every $i\geq0$, the formal group
$\Phi^i_{\mathrm{fl}}(X,\mathbf G_m)$ is smooth, and
$\Phi^i_{\mathrm{\acute et}}(X,\mathbf G_m)$ is prorepresentable,
smooth and isomorphic to it. In particular, this holds for every
abelian variety over $k$.
\end{corollary}

\begin{proof}
For $i\geq1$, combine Proposition~\ref{prop:mazur-ogus-finite-torsion},
Theorem~\ref{thm:artin-mazur-finite-torsion} and
Corollary~\ref{cor:etale-artin-mazur}; for $i=0$, recall that
$G_0=\Phi^0_{\mathrm{\acute et}}(X,\mathbf G_m)$ is smooth. Abelian
varieties have torsion-free crystalline cohomology and an
$E_1$-degenerate Hodge--de Rham spectral sequence.
\end{proof}

The degeneration of the Hodge--de Rham spectral sequence is not needed
for the vanishing of finite torsion: torsion-freeness of crystalline
cohomology suffices, because finite torsion survives in the slope
spectral sequence. For $M\in\Dc$, consider the slope spectral sequence
$E_1^{a,b}=\HH^b(M)^a\Rightarrow\HH^{a+b}(\Tot M)$, and let
$B_\infty^{a,b}\subseteq Z_\infty^{a,b}\subseteq E_1^{a,b}$ be its
infinite boundaries and infinite cycles, so that
$Z_\infty^{a,b}/B_\infty^{a,b}$ is a graded piece of the induced
filtration on $\HH^{a+b}(\Tot M)$. Ekedahl shows that
\[
 B_\infty^{a,b}\subseteq F^\infty B^a\bigl(\HH^b(M)\bigr)
 \subseteq V^{-\infty}Z^a\bigl(\HH^b(M)\bigr)\subseteq Z_\infty^{a,b},
\]
and that every torsion class of
$V^{-\infty}Z^a\bigl(\HH^b(M)\bigr)/B_\infty^{a,b}$ lifts to a torsion
class of $\HH^{a+b}(\Tot M)$ \cite[I, Proposition 2.3]{EkDiag}.

\begin{corollary}\label{cor:torsion-free-crystalline}
Let $n\geq0$. If $\HH^n_{\mathrm{cris}}(X/W)$ is torsion-free, then
$M_X$ has no finite torsion in any bidegree $(a,b)$ with $a+b=n$.
Consequently:
\begin{enumerate}
\item if $\HH^{i+1}_{\mathrm{cris}}(X/W)$ is torsion-free, then
$\Phi^i_{\mathrm{fl}}(X,\mathbf G_m)$ is smooth;
\item if $i\geq1$ and both $\HH^i_{\mathrm{cris}}(X/W)$ and
$\HH^{i+1}_{\mathrm{cris}}(X/W)$ are torsion-free, then
$\Phi^i_{\mathrm{\acute et}}(X,\mathbf G_m)$ is prorepresentable,
smooth and isomorphic to $\Phi^i_{\mathrm{fl}}(X,\mathbf G_m)$.
\end{enumerate}
In particular, if all the crystalline cohomology groups of $X$ are
torsion-free, the conclusions of
Proposition~\ref{prop:mazur-ogus-finite-torsion} and
Corollary~\ref{cor:mazur-ogus-artin-mazur} hold without assuming that
the Hodge--de Rham spectral sequence degenerates.
\end{corollary}

\begin{proof}
By \cite[II, Th\'eor\`eme 1.4]{Illusie},
$\HH^n(\Tot M_X)\simeq\HH^n_{\mathrm{cris}}(X/W)$. Let $a+b=n$, put
$L=\HH^b(M_X)$, and let $x\in V^{-\infty}Z^a(L)$ represent an element
$c$ of $\coeur^{a,b}(X)_{\tors}$. Then $p^mx\in F^\infty B^a(L)$ for
some $m$, and $F^\infty B^a(L)$, being contained in $L_{\tors}$, is
killed by a power of $p$ (see the beginning of
Section~\ref{sec:ordinary}). Hence the class of $x$ in
$V^{-\infty}Z^a(L)/B_\infty^{a,b}$ is torsion, so it is the image of a
torsion class of $\HH^n_{\mathrm{cris}}(X/W)$, which must vanish.
Therefore $x\in B_\infty^{a,b}\subseteq F^\infty B^a(L)$, that is,
$c=0$. This proves the first assertion. Assertion (1) now follows from
Theorem~\ref{thm:artin-mazur-finite-torsion}, assertion (2) from
Corollary~\ref{cor:etale-artin-mazur}, and the last assertion from
\eqref{eq:whole-diagonal}, (1), (2) and the smoothness of $G_0$.
\end{proof}

Assertion (1) recovers, over an arbitrary perfect field, a theorem of
Grammatica \cite[Theorem 3.16]{Grammatica}, who proves it over
algebraically closed fields.

\section{The three dual pairs for a surface}\label{sec:surfaces}

In this section $k$ is an arbitrary perfect field of characteristic
$p>0$, and $S/k$ is a smooth projective geometrically connected
surface. Fix an algebraic closure $\bar k$ of $k$, which is also a
separable closure because $k$ is perfect, and put
\[
 \Gamma_k=\Gal(\bar k/k),\qquad \bar W=W(\bar k),\qquad
 \bar S=S\times_k\bar k.
\]
The results of Sections~\ref{sec:intro}--\ref{sec:artin-mazur} hold
over every perfect field; we use them over $k$ and over $\bar k$. In
particular, Section~\ref{sec:artin-mazur} provides
\[
 M_S\in D^b_c(R),\qquad M_{\bar S}\in D^b_c\bigl(R(\bar k)\bigr),
\]
with $\HH^j(M_S)^i=\HH^j(S,W\Omega_S^i)$, and the modules
$\coeur^{i,j}(S)=\coeur^{i,j}(M_S)$ and
$\coeur^{i,j}(\bar S)=\coeur^{i,j}(M_{\bar S})$. The group $\Gamma_k$
acts on $\bar S$ through the factor $\bar k$, hence by transport of
structure on $\coeur^{i,j}(\bar S)$; this action is semilinear with
respect to $\bar W$ and commutes with $F$ and $V$. For a Dieudonn\'e
module $D$ over $W$, we regard $D\otimes_W\bar W$ as a Dieudonn\'e module
over $\bar W$ with the operators $F\otimes\sigma$ and
$V\otimes\sigma^{-1}$, and we let $\Gamma_k$ act on it semilinearly
through the second factor.

The N\'eron--Severi term below is built from the geometric group
$\NS(\bar S)[p^\infty]$ together with its Galois action, whereas the
two connected Picard terms are Dieudonn\'e modules over $k$ itself.

\subsection{The finite Picard quotient and the two duals}

Let $\Pic^\tau_{S/k}$ be the union of the torsion components of the
Picard scheme, and let $(\Pic^0_{S/k})_{\mathrm{red}}$ be its reduced
identity component, which is an abelian variety because $k$ is
perfect. Let
\begin{equation}\label{eq:finite-picard}
 H=\bigl(\Pic^\tau_{S/k}/(\Pic^0_{S/k})_{\mathrm{red}}\bigr)[p^\infty]
\end{equation}
be the $p$-primary part of the finite quotient. As in
\S\ref{subsec:obstruction}, it decomposes canonically, compatibly with
extension of perfect fields, as
\[
 H=H_{\mathrm{el}}\oplus H_{\mathrm{ll}}\oplus H_{\mathrm{le}},
 \qquad H_{\mathrm{el}}(\bar k)=\NS(\bar S)[p^\infty],
\]
where the second equality is an identification of finite
$\Gamma_k$-modules. The identity component of $H$ is
\[
 H^\circ=H_{\mathrm{ll}}\oplus H_{\mathrm{le}}
 =\Pic^0_{S/k}/(\Pic^0_{S/k})_{\mathrm{red}}.
\]

Put
\[
 \Pi=\NS(\bar S)[p^\infty]=H_{\mathrm{el}}(\bar k),\qquad
 \bar{\mathsf N}=\Pi\otimes_{\ZZ_p}\bar W,
\]
where $\bar{\mathsf N}$ carries the operators $F=1\otimes\sigma$ and
$V=pF^{-1}$, and the semilinear action $g(x\otimes w)=g(x)\otimes g(w)$
of $g\in\Gamma_k$, which commutes with $F$ and $V$. With $\DM_k$ as in
\S\ref{subsec:obstruction}, define the finite-length Dieudonn\'e modules
over $W=W(k)$
\begin{equation}\label{eq:three-modules}
 \mathsf N=\bar{\mathsf N}^{\Gamma_k},\qquad
 \mathsf L=\DM_k(H_{\mathrm{ll}}),\qquad
 \mathsf E=\DM_k(H_{\mathrm{le}}).
\end{equation}
If $p^n\Pi=0$, then $\bar{\mathsf N}$ is a module over $W_n(\bar k)$,
and the action of $\Gamma_k$ on $\Pi$ factors through the Galois group
of a finite Galois extension of $k$. Finite Galois descent therefore
gives a canonical isomorphism, compatible with $F$, $V$ and $\Gamma_k$,
\begin{equation}\label{eq:NS-descent}
 \mathsf N\otimes_W\bar W\simeq\bar{\mathsf N}.
\end{equation}
In general $\mathsf N$ is not $\Pi^{\Gamma_k}\otimes_{\ZZ_p}W$: taking
invariants of $\Pi$ alone loses a nontrivial Galois action. When $k$ is
algebraically closed, $\mathsf N=\NS(S)[p^\infty]\otimes_{\ZZ_p}W$.

Both $F$ and $V$ are nilpotent on $\mathsf L$, while on $\mathsf E$ the
operator $F$ is bijective and $V$ is nilpotent. Since Dieudonn\'e theory
over $k$ is an equivalence of categories, $\mathsf L$ and $\mathsf E$
determine $H_{\mathrm{ll}}$ and $H_{\mathrm{le}}$ as group schemes over
$k$, including their descent data from $\bar k$. With the covariant
convention, $\mathsf N$ is not $\DM_k(H_{\mathrm{el}})$, on which $V$
rather than $F$ is bijective: it is the semi-simple module attached to
the Galois module $\Pi$ of components.

Each dual below has underlying $W$-module $\Hom_W(D,K/W)$, but its
operators depend on the type of $D$. For a semi-simple Dieudonn\'e
module $D$ of finite length, the \emph{unit dual} $D^{\ud}$ is defined
by
\begin{equation}\label{eq:unit-dual}
 (F\lambda)(x)=\sigma\bigl(\lambda(F^{-1}x)\bigr),\qquad V=pF^{-1}.
\end{equation}
For a Dieudonn\'e module $D$ of finite length on which $F$ and $V$ are
nilpotent, the \emph{Dieudonn\'e dual} $D^{\cd}$ is defined by
\begin{equation}\label{eq:dieudonne-dual}
 (F\lambda)(x)=\sigma\bigl(\lambda(Vx)\bigr),\qquad
 (V\lambda)(x)=\sigma^{-1}\bigl(\lambda(Fx)\bigr).
\end{equation}
Both operations commute with extension of perfect fields. In
particular,
\[
 \mathsf N^{\ud}\simeq
 \bigl(\Pi^\vee\otimes_{\ZZ_p}\bar W\bigr)^{\Gamma_k},
\]
where $\Pi^\vee=\Hom_{\ZZ_p}(\Pi,\QQ_p/\ZZ_p)$ carries the
contragredient action $(g\lambda)(x)=\lambda(g^{-1}x)$ and $F$ is
induced by $1\otimes\sigma$. The two superscripts are essential:
\eqref{eq:unit-dual} and \eqref{eq:dieudonne-dual} are different
duality operations.

\subsection{The classification}

\begin{theorem}[Finite torsion of a surface over a perfect field]\label{thm:surface}
Let $\mathsf N$, $\mathsf L$ and $\mathsf E$ be as in
\eqref{eq:three-modules}. The only possibly nonzero finite-torsion
pieces are the following three dual pairs of Dieudonn\'e modules over
$W(k)$, where $(i,j)$ denotes the position $\HH^j(S,W\Omega_S^i)$:
\begin{center}
\upshape
\renewcommand{\arraystretch}{1.25}
\begin{tabular}{@{}llll@{}}
\toprule
Part & First position and module & Dual position and module & Origin\\
\midrule
semi-simple & $(1,1)$:\ $\mathsf N$ & $(1,2)$:\ $\mathsf N^{\ud}$
 & $\NS(\bar S)[p^\infty]$\\
nilpotent & $(0,2)$:\ $\mathsf L$ & $(1,2)$:\ $\mathsf L^{\cd}$
 & $H_{\mathrm{ll}}$\\
semi-simple & $(0,2)$:\ $\mathsf E$ & $(2,1)$:\ $\mathsf E^{\ud}$
 & $H_{\mathrm{le}}$\\
\bottomrule
\end{tabular}
\end{center}
All pieces not listed are zero. The same statements hold for
$A_f^2$ and $A_{\SSpart}^1$, for the nilpotent and semi-simple parts of
c\oe ur torsion, and for the modules
$\HH^{-i}\bigl(t^1\Hd^{i+j}(M_S)\bigr)^i$ of
Theorem~\ref{thm:comparison}.
\end{theorem}

The positions of the three pairs overlap. Equivalently, the complete
table of c\oe ur torsion is
\begin{equation}\label{eq:surface-table}
\begin{array}{c|ccc}
 \coeur^{i,j}(S)_{\tors} & i=0 & i=1 & i=2\\ \hline
 j=2 & \mathsf L\oplus\mathsf E & \mathsf L^{\cd}\oplus\mathsf N^{\ud} & 0\\[2pt]
 j=1 & 0 & \mathsf N & \mathsf E^{\ud}\\[2pt]
 j=0 & 0 & 0 & 0.
\end{array}
\end{equation}

\begin{proof}
We first compare semi-simple parts over $k$ and over $\bar k$, then
identify the three initial modules, and finally apply duality over $k$.

\emph{Base change and descent.}
Only the semi-simple parts need to be compared with their analogues over
$\bar k$. By Theorem~\ref{thm:comparison} and \eqref{eq:A-pieces},
$\coeur^{i,j}(M)_{\SSpart}\simeq A^1\bigl(\HH^j(M)\bigr)^i_{\tors}$ for
$M=M_S$ and for $M=M_{\bar S}$. Ekedahl's completed pullback along
$k\to\bar k$, from $D^b_c(R)$ to $D^b_c(R(\bar k))$, is exact and
conservative on coherent modules \cite[III, Lemma 1.4]{EkII}; it carries
$M_S$ to $M_{\bar S}$, and it commutes with internal derived Hom and with
the duality functor $\DD$ \cite[I, \S9]{EkII}. It therefore carries the
biduality spectral sequence of $\HH^j(M_S)$ to that of $\HH^j(M_{\bar
S})$, and hence $A^1(\HH^j(M_S))$ to $A^1(\HH^j(M_{\bar S}))$. The
modules $A^1(\HH^j(M))^i$ are finitely generated over $W$
\cite[IV, Theorem 3.3]{EkI}, and on such modules completed pullback is
ordinary extension of scalars, which preserves torsion submodules.
Hence there are natural isomorphisms
\begin{equation}\label{eq:finite-torsion-base-change}
 \coeur^{i,j}(S)_{\SSpart}\otimes_W\bar W
 \simeq\coeur^{i,j}(\bar S)_{\SSpart},
\end{equation}
compatible with $F$, $V$ and the semilinear actions of $\Gamma_k$. No
base change statement is needed for the whole Hodge--Witt modules,
which need not be finitely generated over $W$. As
$(D\otimes_W\bar W)^{\Gamma_k}=D$ for every $W$-module $D$ of finite
length, it follows that
$\coeur^{i,j}(S)_{\SSpart}\simeq(\coeur^{i,j}(\bar S)_{\SSpart})^{\Gamma_k}$.

\emph{The N\'eron--Severi module.}
We first work over $\bar k$. For every $n\geq1$, the Kummer sequence
gives an exact sequence
\[
 0\longrightarrow\Pic(\bar S)/p^n
 \longrightarrow\HH^2_{\mathrm{fl}}(\bar S,\mu_{p^n})
 \longrightarrow\Br(\bar S)[p^n]\longrightarrow 0.
\]
The group $\Pic^0(\bar S)$ is the group of $\bar k$-points of the
abelian variety $(\Pic^0_{\bar S/\bar k})_{\mathrm{red}}$, so it is
divisible and $\Pic(\bar S)/p^n=\NS(\bar S)/p^n$. The transition maps
of this system are surjective, so passing to the limit gives an exact
sequence
\begin{equation}\label{eq:kummer-NS}
 0\longrightarrow\NS(\bar S)\otimes_{\ZZ}\ZZ_p
 \longrightarrow\HH^2_{\mathrm{fl}}(\bar S,\ZZ_p(1))
 \longrightarrow T_p\Br(\bar S)\longrightarrow 0,
\end{equation}
where
$\HH^2_{\mathrm{fl}}(\bar S,\ZZ_p(1))=\varprojlim_n\HH^2_{\mathrm{fl}}(\bar S,\mu_{p^n})$
and $T_p\Br(\bar S)=\varprojlim_n\Br(\bar S)[p^n]$, the transition maps
of the last system being multiplication by $p$. A Tate module is
torsion-free: if $p^r(x_n)_n=0$, then $x_n=p^rx_{n+r}=0$ for every $n$.
Since $\NS(\bar S)$ is finitely generated, \eqref{eq:kummer-NS} gives
\begin{equation}\label{eq:kummer-torsion}
 \HH^2_{\mathrm{fl}}(\bar S,\ZZ_p(1))_{\tors}
 \simeq\NS(\bar S)[p^\infty]=\Pi.
\end{equation}

The map $d\log$ induces an isomorphism
$\HH^2_{\mathrm{fl}}(\bar S,\ZZ_p(1))\simeq
\HH^1(\bar S,W\Omega^1_{\bar S,\log})$ \cite[I, Corollary 3.27]{Illusie}.
By the theorem of Illusie and Raynaud on logarithmic cohomology
\cite[IV, Theorem 3.3]{IR}, the natural map
\[
 \HH^1(\bar S,W\Omega^1_{\bar S,\log})
 \longrightarrow\coeur^{1,1}(\bar S)^{F=1}
\]
is surjective, and its kernel is the group of $\bar k$-points of a
connected unipotent group of dimension $T^{0,1}$, the domino number of
the domino
$\HH^1(\bar S,W\mathcal O_{\bar S})/V^{-\infty}Z^0\to F^\infty B^1$
arriving in bidegree $(1,1)$. Ekedahl's duality for domino numbers
\cite[IV, Corollary 3.5.1]{EkI} gives $T^{0,1}=T^{0,3}$, and $T^{0,3}=0$
because $\HH^3(M_{\bar S})=0$. Hence the map above is an isomorphism.
Taking torsion subgroups and using \eqref{eq:kummer-torsion}, we obtain
\[
 \Pi\simeq\bigl(\coeur^{1,1}(\bar S)_{\tors}\bigr)^{F=1}
 =\bigl(\coeur^{1,1}(\bar S)_{\SSpart}\bigr)^{F=1},
\]
because the nilpotent part has no nonzero $F$-fixed elements. For a
semi-simple Dieudonn\'e module $D$ of finite length over $\bar W$, the
natural map $D^{F=1}\otimes_{\ZZ_p}\bar W\to D$ is an isomorphism; in
the language of $p$-adic Hodge theory, \'etale $\varphi$-modules of
finite length over $W(\bar k)$ are trivial. If $pD=0$, this is
\cite[Proposition 1.1]{KatzSGA7}, whose proof rests on Lang's theorem
for $\mathrm{GL}_n$. In general one argues by induction on the length,
using the exact sequence $0\to pD\to D\to D/pD\to0$: by
\cite[Proposition 1.2]{KatzSGA7} and induction, $1-F$ is surjective on
$D/pD$ and on $pD$, hence on $D$, so the sequence remains exact on
$F$-invariants and the five lemma applies. Therefore
\[
 \coeur^{1,1}(\bar S)_{\SSpart}\simeq\Pi\otimes_{\ZZ_p}\bar W
 =\bar{\mathsf N}.
\]
All the maps involved are natural with respect to automorphisms of
$\bar k$ over $k$, so this isomorphism is $\Gamma_k$-equivariant as well
as compatible with $F$ and $V$. Taking $\Gamma_k$-invariants and using
the base change step, we obtain
\begin{equation}\label{eq:NS-identification}
 \coeur^{1,1}(S)_{\SSpart}\simeq\bar{\mathsf N}^{\Gamma_k}=\mathsf N.
\end{equation}

\emph{The connected finite Picard module.}
By \eqref{eq:H1-WO} and Corollary~\ref{cor:etale-artin-mazur}, the
functor $\Phi^1_{\mathrm{\acute et}}(S,\mathbf G_m)$ is prorepresentable
and equal to $G_1=\Phi^1_{\mathrm{fl}}(S,\mathbf G_m)$; it is the formal
Picard group of $S$, the formal completion of $\Pic_{S/k}$ at the
origin \cite[II, \S4]{AM}. Its reduced subgroup is the formal completion
of the abelian variety $(\Pic^0_{S/k})_{\mathrm{red}}$, and the quotient
is the infinitesimal group scheme
$\Pic^0_{S/k}/(\Pic^0_{S/k})_{\mathrm{red}}$; thus $Q_1=H^\circ$.
Theorem~\ref{thm:artin-mazur-finite-torsion} with $X=S$ and $i=1$
therefore gives, over $k$ itself,
\begin{equation}\label{eq:connected-picard-seeds}
 \coeur^{0,2}(S)_{\tors}=\HH^2(S,W\mathcal O_S)[V^\infty]
 \simeq\DM_k(H^\circ)=\mathsf L\oplus\mathsf E,
\end{equation}
with $\coeur^{0,2}(S)_{\nil}\simeq\mathsf L$ and
$\coeur^{0,2}(S)_{\SSpart}\simeq\mathsf E$.

\emph{Duality and vanishing.}
For a smooth proper geometrically connected variety $X/k$ of dimension
$N$, Ekedahl's geometric duality over the perfect field $k$
\cite[IV, Theorem 3.5]{EkI} (see also \cite[II, Theorem 1.1(ii)]{EkII}),
combined with the formulas of \cite[IV, Theorem 3.3]{EkI} and
Theorem~\ref{thm:comparison}, reads
\begin{align}
 \coeur^{i,j}(X)_{\SSpart}&\simeq
 \bigl(\coeur^{N-i,N+1-j}(X)_{\SSpart}\bigr)^{\ud},
 \label{eq:geometric-unit-duality}\\
 \coeur^{i,j}(X)_{\nil}&\simeq
 \bigl(\coeur^{N-1-i,N+2-j}(X)_{\nil}\bigr)^{\cd};
 \label{eq:geometric-nil-duality}
\end{align}
the internal shifts in Ekedahl's formulas are already built into these
indices, and the duals carry the operators \eqref{eq:unit-dual} and
\eqref{eq:dieudonne-dual}. For $N=2$, the first formula pairs $(1,1)$
with $(1,2)$ and $(0,2)$ with $(2,1)$, and the second pairs $(0,2)$
with $(1,2)$. Together with \eqref{eq:NS-identification} and
\eqref{eq:connected-picard-seeds}, this gives all the modules in the
table.

It remains to see that no other position contributes. Cohomology
vanishes outside $0\leq i,j\leq2$. The groups $\HH^0(S,W\Omega_S^i)$
are $p$-torsion-free, because the sheaves $W\Omega_S^i$ are
\cite[I, \S3]{Illusie}; by Proposition~\ref{prop:ordinary-comparison},
their c\oe urs are torsion-free. Next, $\coeur^{0,1}(S)_{\tors}=0$. This
is \eqref{eq:H1-WO} for $X=S$, but it can also be seen directly over $k$.
Since $S$ is proper and geometrically connected, the reduction map
$\HH^0(S,W\mathcal O_S)=W(k)\to\HH^0(S,\mathcal O_S)=k$ is surjective,
so the exact sequence
\[
 0\longrightarrow W\mathcal O_S\xrightarrow{\ V\ }W\mathcal O_S
 \longrightarrow\mathcal O_S\longrightarrow 0
\]
shows that $V$ is injective on $\HH^1(S,W\mathcal O_S)$. As nothing
arrives in internal degree zero under $d$, the module
$\coeur^{0,1}(S)_{\tors}$ is a submodule of $\HH^1(S,W\mathcal O_S)$; it
has finite length and $V$ is topologically nilpotent on it, so it is
killed by a power of $V$ and therefore vanishes. For the nilpotent part,
\eqref{eq:geometric-nil-duality} makes every position with $j=1$ dual
to one with $j=3$, and every position with $i=2$ dual to one with
$i=-1$; all of these vanish, leaving only $(0,2)$ and $(1,2)$. For the
semi-simple part, the vanishing at $(0,1)$ gives vanishing at $(2,2)$ by
\eqref{eq:geometric-unit-duality}, leaving exactly the four positions
listed.
\end{proof}

\begin{corollary}\label{cor:surface-diagonal}
For a surface $S$ as above, $t^1\Hd^n(M_S)=0$ unless $n=2$ or $3$, and
\begin{align*}
 t^1\Hd^2(M_S)&\simeq
 (\mathsf L\oplus\mathsf E)\oplus\mathsf N(-1)[1],\\
 t^1\Hd^3(M_S)&\simeq
 (\mathsf L^{\cd}\oplus\mathsf N^{\ud})(-1)[1]
 \oplus\mathsf E^{\ud}(-2)[2].
\end{align*}
Thus the finite-torsion contributions in total degrees $2$ and $3$ are
the $W$-modules
\[
 \mathsf L\oplus\mathsf E\oplus\mathsf N
 \qquad\text{and}\qquad
 \mathsf L^{\cd}\oplus\mathsf N^{\ud}\oplus\mathsf E^{\ud}.
\]
These statements concern diagonal finite torsion only; they say nothing
about the torsion carried by dominos.
\end{corollary}

\begin{proof}
Insert \eqref{eq:surface-table} into
Corollary~\ref{cor:whole-diagonal}. The only nonzero entries have
$i+j=2$ or $3$.
\end{proof}

\begin{remark}\label{rem:formal-brauer-surface}
By Corollary~\ref{cor:etale-artin-mazur} and
\eqref{eq:connected-picard-seeds}, the formal Brauer group
$\Phi^2_{\mathrm{\acute et}}(S,\mathbf G_m)$ of the surface $S$ is
prorepresentable if and only if $\mathsf L=\mathsf E=0$, that is, if
and only if $H^\circ=0$, i.e.\ $\Pic_{S/k}$ is reduced. It is then also
smooth, since $\HH^3(S,W\mathcal O_S)=0$.
\end{remark}

\section*{Declarations}

\paragraph{Relation to the author's thesis.}
Part of this work appeared, over algebraically closed fields, in the
author's PhD thesis~\cite{YangThesis}; the present paper gives complete
proofs and works over an arbitrary perfect field.

\paragraph{Funding.}
Y.Y.~is partially supported by a grant from the National Natural Science
Foundation of China (NSFC-12231001).

\paragraph{Conflict of interest.}
The author declares that there is no conflict of interest.

\paragraph{Use of AI tools.}
To turn part of the author's thesis and proofs into a self-contained
article, ChatGPT~6~Astra (OpenAI) was used to generate a first draft of
this paper from this material. Claude Code (Anthropic), with the Claude
Opus~5.5 model, was then used to revise the draft into the present
manuscript, including the exposition and the written form of the
proofs. The author has checked all arguments and takes full
responsibility for the content of this paper.

\end{document}